\documentclass{article}
\usepackage{graphicx} % Required for inserting images
\usepackage[top=2cm, bottom=2cm, right=2cm, left=2cm]{geometry}
\usepackage{caption}
\usepackage{subcaption}
\usepackage{algorithm}
\usepackage{algpseudocode}
\usepackage{authblk}
\usepackage{orcidlink} 
\usepackage{amsfonts}
\usepackage{amsmath}
\usepackage{float}

\usepackage{amsthm}
\theoremstyle{remark}
\newtheorem{remark}{Remark} 

\theoremstyle{acknowledgement}
\newtheorem*{acknowledgement}{Acknowledgement}

\title{\textbf{Parameter selection procedures for model reduction of a parameterized wave propagation in damaged material}}
\author[1]{\normalsize{\textbf{Hannah Burmester}}
\orcidlink{0009-0003-9802-6902}
\footnote{ e-mail \texttt{hannah.burmester@tu-braunschweig.de},
            phone +49\,531\,391\,7428.
            }}
\author[1]{\textbf{Carmen Grä\ss{}le} \orcidlink{0000-0003-0318-0740}
\footnote{e-mail \texttt{c.graessle@tu-braunschweig.de},
            phone +49\,531\,391\,7417.
            }}
\author[1]{\textbf{Saddam Hijazi} \orcidlink{0000-0002-9114-1102}
\footnote{e-mail \texttt{saddam.hijazi@tu-braunschweig.de},
            phone +49\,531\,391\,7484.
            }}
\affil[1]{%\CountryCode[DE]
Technische Universität Braunschweig, Institute for Partial
Differential Equations, Universitätsplatz 2, 38106 Braunschweig}

\date{October 2026}

\providecommand{\WileyBibTextsc}{}
\let\textsc\WileyBibTextsc
\providecommand{\othercit}{}
\providecommand{\jr}[1]{#1}

\begin{document}
\maketitle
\begin{abstract}

In this work, we consider a parameterized version of the classical wave equation that models wave propagation in damaged material. 
We apply model order reduction to the wave equation using a reduced basis method with the POD-Greedy algorithm. 
The success of the model order reduction highly depends on the parameter set used in the POD-Greedy algorithm. While a larger parameter set leads to more insight on the behavior of the full order model, it also prolongs the reduction process in the offline phase for constructing the reduced basis. Therefore, subsampling strategies based on Discrete Empirical Interpolation Method, pivoted QR decomposition and Polytope Division Method are applied. 
\end{abstract}
%% maketitle must follow the abstract.
%\maketitle                   % Produces the title.

\section{Introduction}

Parameter dependent partial differential equations (PPDE) are used to model various real-world phenomena. Techniques such as finite volumes, finite elements and finite differences are used to solve the PPDE numerically. Solving the resulting high dimensional systems of equations often leads to high computational efforts, that can be challenging for example in multi-query scenarios where the PPDE has to be solved repeatedly for different parameter inputs. The high computational effort motivates the use of model order reduction techniques. The basic procedure in model order reduction is to replace a high dimensional model with a model of reduced dimension while preventing too much loss in accuracy.\\
Here, we consider a classical wave equation modeling wave propagation in damaged material. A damage is modeled as a local reduction in the wave propagation velocity, causing the reflection of a propagating wave at the damage. The wave equation is discretized in space and time using finite elements and the $\theta$-method, respectively. At this point, we would like to contextualize our problem within the framework of damage identification, where we aim to solve an inverse problem to determine the unknown parameters modeling the damage by utilizing both the forward model and sensor measurements. Due to the fact that such a process inherently requires several forward model evaluations, leveraging reduced order models becomes essential to mitigate the overall computational costs, see e.g.\ \cite{MasterthesisHannah,GH2024,HFBRG25,Abschlussband2025,MGRMLL23}.\\
For the model order reduction of the high dimensional model, a reduced basis method with a Greedy approach using Proper Orthogonal Decomposition (POD) is implemented. The POD-Greedy procedure\cite{HO2008,RBPDE,Haas17,MP21} is used to create a reduced basis (RB) in an iterative manner. Based on the reduced basis and the high dimensional model, the reduced model is constructed using the RB-Galerkin method. We refer to \cite{Glas2020} for the reduced basis method for the wave equation. \\
Model order reduction is challenging for hyperbolic equations such as the considered wave equation due to the limitation associated with the Kolmogorov barrier \cite{transportdominatedMOR, RBlimitation}. Moreover, the success of the model order reduction highly depends on choosing a suitable parameter set during the reduction process in the offline phase of the reduced basis generation. While a larger parameter set (i.e.\ a finer discretization of the parameter domain) leads to more insight on the behavior of the full order model, it also prolongs the reduction process. Therefore, strategies are investigated to simplify the model order reduction process with suitable parameter subsampling methods. We apply subsampling strategies based on the Discrete Empirical Interpolation Method (DEIM) presented in \cite{Subsampling}, the pivoted QR decomposition introduced in \cite{Subsampling} and strategies based on the Polytope Division Method described in \cite{nielen2024polytope}.

\section{The wave model}
\label{sec:waveModel}

We consider the classical wave equation, in which damage is incorporated via parameterization. The current study case is inspired by the model used in \cite{BML2021} where a damage is characterized by its location and a local reduction of the wave propagation velocity. Here, an extension of this model is studied which has been proposed in \cite{MasterthesisHannah}: an additional parameter related to the size of the damage is introduced and  the spatial discretization is changed to the finite element method requiring a smoothing modification of the initial condition and the damage parameterization. We further note that this model is reported in \cite{GH2024} and some selected examples for POD-Greedy model reduction are presented. Further, the model is studied in \cite{HFBRG25} as a test case for a comparison of intrusive and non-intrusive model reduction for second order systems. \\
\\
Let $\Omega \subset \mathbb{R}^2$ be the spatial domain with sufficiently smooth boundary $\partial \Omega$ and let $[0,T] = [0,5]$ be the time interval of interest. The time interval is separated into two phases. In the excitation phase, $t\in [0,0.1],$ the displacement $u(t,\textbf{x})$ at the center $\textbf{x}_C = (2.5,2.5)\in \Omega$ is given as 
\begin{equation*}
     u_{\textbf{x}_C}(t) = 10 + \sin (2\pi At), \quad A = 450,000
\end{equation*}
and smoothed outwards via multiplication with 
\begin{equation*}
    g(\textbf{x}) = 
    \Bigg\{
    \begin{array}{cl}
            \exp \left( 1+ \frac{0.09}{(x_1-2.5)^2+(x_2-2.5)^2 -0.09}  \right), &  \| \textbf{x}-\textbf{x}_C \|_2 < 0.3, \\
             0, & \| \textbf{x}- \textbf{x}_C \|_2 \geq 0.3,
    \end{array}
\end{equation*}
where $\textbf{x} = \left( x_1, x_2 \right)$.
After the excitation phase, for $t\in [0.1,5.0],$ the wave propagation is modeled with the initial boundary value problem 
\begin{equation}\label{eq:ibvp_wave}
\left\{\begin{aligned}
 u_{tt}(t,\textbf{x}) - c_\mu^2(\textbf{x}) \Delta u(t,\textbf{x}) &= 0, && \quad \text{in } (0.1,5.0] \times \Omega,  \\
  u(t,\textbf{x}) & = 0, && \quad \text{on } (0.1,5.0] \times \partial\Omega,\\
 u(0.1,\textbf{x}) &= u_0(\textbf{x}) , &&\quad \text{in } \Omega, \\
 \dot{u}(0.1,\textbf{x}) &= u_1(\textbf{x}) , && \quad \text{in } \Omega,
\end{aligned}\right.
\end{equation}
with displacement $u(t,\textbf{x})$, parameter dependent wave propagation velocity $c_\mu(\textbf{x})$ and initial conditions
\begin{align*}
    u_0 (\textbf{x}) &= u_{\textbf{x}_C} (0.1) \cdot g(\textbf{x})  \qquad \text{ and }\\
        u_1 (\textbf{x}) &= \dot{u}_{\textbf{x}_C} (0.1) \cdot g(\textbf{x}).
\end{align*}
The parameter dependent wave propagation velocity is given by
\begin{equation*}
    c^2_\mu (\textbf{x}) = 0.25-(0.25-c_\text{D}) \exp \left( - 
        \frac{(x_1-x_\text{1D})^2 + (x_2-x_\text{2D})^2}{d_\text{D}}\right), 
\end{equation*}
such that the wave propagation velocity at the center of the damage $(x_\text{1D},x_\text{2D})\in \Omega$ is reduced to $c_\text{D}$ in comparison to $c^2 = 0.25$ in the intact material. For this work, we assume a fixed damage location $(x_\text{1D},x_\text{2D})$. The parameter $d_\text{D}$ influences the size of the damage. The parameters characterizing the damage are combined to the parameter vector 
\begin{equation*}
    \mu = (c_\text{D}, d_\text{D}) \in \mathcal{P}
\end{equation*}
with parameter domain 
\begin{equation*}
    \mathcal{P} = [0.02,0.20] \times [0.05,0.50].
\end{equation*} 
We multiply \eqref{eq:ibvp_wave} with a test function $v\in V = H_0^1 (\Omega),$ integrate over space and arrive at the weak formulation 
\begin{equation}\label{eq:ibvp_wave:weak}
\left\{\begin{aligned}
     \langle \ddot{u}(t),v \rangle_{V' \times V} + (c_\mu^2 (\textbf{x}) \nabla u(t), \nabla v )_{L_2 (\Omega )} &= 0   && \quad \forall v \in V, t \in  (0.1,5.0]  \; a.e., \,  \\
 u(0.1) &= u_0 , &&\\
 \dot{u}(0.1) &= u_1. && 
\end{aligned}\right. 
\end{equation}

\begin{figure}
    \centering
    \includegraphics[width=0.75\linewidth]{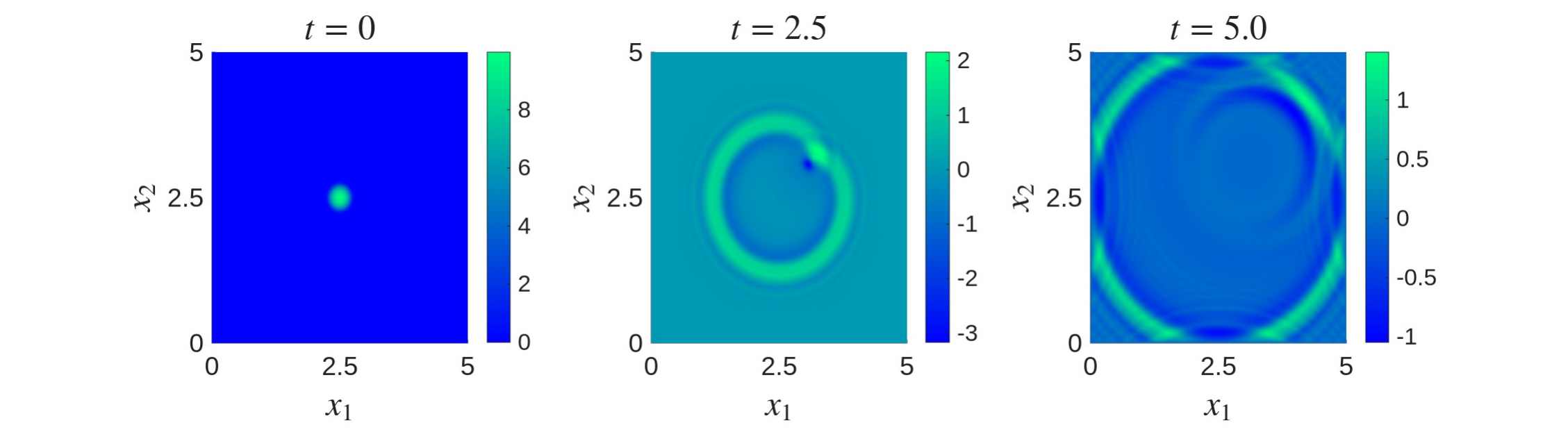}
    \caption{High dimensional solution of the wave equation with damage location at $(x_\text{1D},x_\text{2d})=(3,3)$ and $\mu = (c_\text{D}, d_\text{D}) = (0.05,0.1)$ at times $t = 0,2.5,5.0$ (f.l.t.r.).}
    \label{fig:wellenmodell:hdimLsg}
\end{figure}

To discretize the weak formulation \eqref{eq:ibvp_wave:weak} in space, we use piecewise linear finite elements with $\psi_i(\textbf{n}_j) = \delta_{ij}$ where $\{\textbf{n}_1,..., \textbf{n}_N\}$ denotes the set of nodes of the finite element mesh.
We search for an approximation 
\begin{equation*}
    u_h(t,\textbf{x}) = \sum_{i=1}^{N} (\textbf{u}_h (t))_i \psi_i(\textbf{x}) \approx u(t,\textbf{x})
\end{equation*}
in the finite element space $V_h = \text{span}\langle \psi_1,...,\psi_N\rangle \subset V$
and restrict the test space in \eqref{eq:ibvp_wave:weak} to $V_h.$ This leads to 
\begin{equation}\label{eq:ibvp_wave:discretized_in_space}
     \left\langle \sum_{i=1}^{N} (\ddot{\textbf{u}}_h(t))_i\psi_i, \psi_j  \right\rangle_{V'\times V}
    + \left\langle c_\mu^2(\textbf{x}) \sum_{i=1}^{N} (\textbf{u}_h(t))_i \nabla \psi_i, \nabla \psi_j \right\rangle_{L_2(\Omega)} 
    = 0,   
       \quad j=1,\dots, N, t \in  (0.1,5.0] \; a.e., 
\end{equation}
together with the initial conditions $u_h(0.1) = u_{0,h}$ and $\dot{u}_h(0.1) =  u_{1,h}$, where $u_{0,h}$ and $u_{1,h}$ denote suitable finite element approximations of $u_0$ and $u_1$, respectively. From \eqref{eq:ibvp_wave:discretized_in_space} we deduce the system of ordinary differential equations (ODEs)
\begin{equation}\label{eq:wave:FOM:ODE}
M \ddot{\textbf{u}}_h(t) + S \textbf{u}_h(t) = 0,   
       \qquad  t \in  (0.1,5.0], 
\end{equation}
with $M_{i,j} := \langle \psi_i , \psi_j \rangle_{L_2(\Omega )}$ and $S_{i,j} := \langle c_\mu^2(\textbf{x}) \nabla \psi_i, \nabla \psi_j \rangle_{L_2 (\Omega)}$ denoting the mass matrix and the stiffness matrix, respectively. 
In the numerical tests, we use a finite element mesh consisting of $N = 11,282$ degrees of freedom. \\
\\
The resulting system of ODEs \eqref{eq:wave:FOM:ODE} is integrated using a $\theta$-scheme with $\theta=\frac{1}{4}$ following \cite{Glas2020}:

\begin{align*}
    \left( M + \frac{1}{4} \tau^2 S \right) \textbf{u}_h^{k+1} = \,
    &M \left( 2 \textbf{u}_h^k - \textbf{u}_h^{k-1} \right) +
    \tau^2 S \left( -\frac{1}{2} \textbf{u}_h^k - \frac{1}{4} \textbf{u}_h^{k-1}\right).
\label{eq:GlSys:hdimModell}
\end{align*}
For the numerical runs, we use $m = 250$ time points $t_i, i = 1,...,m$, with equidistant step size $\tau=0.02$.  The vector $\textbf{u}_h^k$ approximates $\textbf{u}_h (t_k).$ In the following, we refer to the finite element approximation of the wave equation as the full order model (FOM) or high dimensional model.  \\
Figure~\ref{fig:wellenmodell:hdimLsg} shows the high dimensional solution $u_h(t,\textbf{x})$ at $t = 0,2.5,5.0$ (f.l.t.r.) for a damge located at $(x_\text{1D},x_\text{2D})=(3,3)$ and $\mu = (0.05,0.1).$
We see that the excited wave propagates outward from the center and surpasses the damage at around $t = 2.5$ leading to a second overlaid wave. 

\section{Model order reduction}

The main idea in model order reduction (MOR) is to replace a high dimensional model, also referred to as full order model (FOM), with a reduced order model (ROM), that requires less computational time and less storage capacities than the full order model. The ROM is set up in the so-called \textit{offline phase} which is computationally demanding and involves simulating the high dimensional model. In the following \textit{online phase}, the ROM is used to compute approximations of the high dimensional solutions quickly and with low computational effort. \\
Reduced basis (RB) methods \cite{HO2008,RBPDE,Haas17,MP21} are commonly used in the field of model order reduction, being effectively applicable when the space spanned by the set of all high-dimensional solutions
\begin{equation*}
\mathcal{M} = \text{span} \{u_h(t;\mu): \mu \in \mathcal{P}\}
\end{equation*}
can be well approximated by a subspace of far smaller dimension. Setting up the reduced model consists of two steps: first a reduced basis space of low dimension is determined that well approximates $\mathcal{M}$. Second, the ROM is constructed based on the RB space and the high dimensional model. For step one, we use the POD-Greedy algorithm explained in Section \ref{subsec:PODGreedy} and for step two the RB-Galerkin method described in Section \ref{subsec:GlaerkinRB}. 

\begin{remark}
It is well known that MOR is challenging for hyperbolic problems such as the wave equation, due to the slow decay of the Kolmogorov n-width, see e.g.\ \cite{transportdominatedMOR, RBlimitation}.\\
Moreover, we underline the fact that the wave equation model could be reformulated as first order system with a Hamiltonian structure, where the Hamiltonian represents the conservation of the total physical energy. Standard MOR methods often result in reduced order systems that do not preserve the physical quantities that remain invariant during the evolution of Hamiltonian systems. In particular, standard orthogonal Galerkin projection does not generally preserve the symplectic structure of Hamiltonian systems, which may lead to unstable reduced order models, especially for long-time integration. Therefore, the stability and structure-preserving properties of the reduced order model should be verified a posteriori by evaluating the conserved energy of the system and ensuring that it remains constant over time. In addition, we would like to mention several approaches that have been developed to address the instabilities associated with standard MOR techniques when applied to Hamiltonian systems. These methods construct symplectic reduced order bases that preserve the underlying Hamiltonian structure at the reduced order level \cite{lall2003structure, peng2016symplectic,farhat2015structure, afkham2017structure}. A data-based approach for preserving the incremental stabilty also in nonlinear systems has been proposed in \cite{Abbasi2022}.
\end{remark}

\subsection{RB-Galerkin method}
\label{subsec:GlaerkinRB}
Let us assume to be given a set of orthonormal basis functions that span a low dimensional subspace $V_r=\text{span}\langle \varphi_1,\dots,\varphi_{N_r}\rangle \subset V_h$ of the full order space $V_h$. For the construction of such a basis we refer to Section~\ref{subsec:PODGreedy}.
The reduced basis space $V_r$ substitutes $V_h$ in \eqref{eq:ibvp_wave:discretized_in_space}. We search for an approximation
\begin{equation*}
    u_a(t,\textbf{x}) := \sum_{i=1}^{N_r} (\textbf{u}_r(t))_i \varphi_i(\textbf{x})
\end{equation*}
of $u(t,\textbf{x})$ in $V_r$ and test with the basis functions $\{\varphi_i\}_{i=1}^{N_r}$ of $V_r$ \cite{HO2008,RBPDE,Haas17,MP21}. This leads to the reduced order model (ROM):  
\begin{equation}\label{eq:ibvp_wave:ROM}
\left\{\begin{aligned}
\left\langle \ddot{u}_a(t), \varphi_j \right\rangle_{V'\times V}
    &+ \left\langle c_\mu^2(\textbf{x}) \nabla u_a(t), \nabla \varphi_j \right\rangle_{L_2(\Omega)} 
    = 0,   
       && j = 1,..., N_r, \; t \in  (0.1,5.0] \; a.e., \\
 u_a(0.1) &= u_{0,r}, && \text{in } \Omega, \\
 \dot{u}_a(0.1) &= u_{1,r},  && \text{in } \Omega.
\end{aligned}\right.
\end{equation}
Here $u_{0,r}$ and $u_{1,r}$ are suitable approximations of $u_0 (\textbf{x})$ and $u_1(\textbf{x})$ in the reduced space $V_r$.
We deduce the reduced system of ODEs
\begin{equation}\label{eq:wave:ROM:ODE}
M_r \ddot{\textbf{u}}_r(t) + S_r \textbf{u}_r(t) = 0,   
       \quad t \in  (0.1,5.0), 
\end{equation}
with the reduced matrices $(M_r)_{i,j} := \langle \varphi_i,\varphi_j\rangle_{L_2(\Omega)}$ and $(S_r)_{i,j} := \langle c_\mu^2(\textbf{x}) \nabla \varphi_i,\nabla \varphi_j\rangle_{L_2(\Omega)}$. Since the number of equations in \eqref{eq:wave:ROM:ODE} is reduced to $N_r \ll N$, 
the system \eqref{eq:wave:ROM:ODE} can generally be solved way faster than system \eqref{eq:wave:FOM:ODE}.

To integrate the resulting system of ordinary differential equations with respect to time, the same $\theta$-method is applied using the same time discretization with $\tau = 0.02, T= 5, m = 250$ as for the full order model.

\subsection{The POD-Greedy algorithm }
\label{subsec:PODGreedy}

The POD-Greedy algorithm \ref{Alg:PODGreedy} combines a classical greedy procedure with Proper Orthogonal Decomposition (POD). The greedy algorithm iteratively identifies parameter values for which new basis vectors are determined via POD as described in the following, compare also \cite{HO2008,RBPDE,Haas17,MP21}. \\
\\
First, a finite parameter subset $\mathcal{P}_\text{train} \subset \mathcal{P}$, a maximum number of iterations $N_\text{max}$, an error indicator $\Delta_e(V_N,\mu)$ and an error tolerance $\varepsilon$ are chosen. The POD-Greedy algorithm terminates when the error tolerance and/or the maximum number of iterations $N_\text{max}$ is reached. Initially, a random parameter value $\mu_1 \in \mathcal{P}_\text{train}$ is picked and the corresponding high dimensional solution 
\begin{equation*}
U_h (\mu_1) := [\textbf{u}_h^1(\mu_1), \dots,\textbf{u}_h^m(\mu_1)]    
\end{equation*}
is computed. Based on this solution $U_h (\mu_1)$, $l$ orthonormal basis vectors $\{\varphi_1^i \}_{i=1}^l = : \Phi_1$ are determined via POD. Further basis elements are added iteratively: each iteration begins by determining the next parameter value $\mu_{N+1}\in \mathcal{P}_\text{train}$ that maximizes the considered error indicator $\Delta_e(V_N,\mu)$ on $\mathcal{P}_\text{train}.$ Next, the high dimensional model is solved for the determined parameter value $\mu_{N+1}$. Then, the orthogonal projection of the solution $U_h (\mu_{N+1})$ onto $V_N = \text{span} \langle \Phi_N\rangle$ is subtracted from the solution itself and POD is applied yielding $l$ POD modes $\{ \varphi_{N+1}^i\}_{i=1}^l$. The POD modes $\{ \varphi_{N+1}^i\}_{i=1}^l$ are orthonormalized against the current basis $\Phi_N$ and added to the current reduced basis.

\begin{algorithm}
	\caption{POD-Greedy, compare \cite{ HO2008,Haas17,Glas2020} } 
    \label{Alg:PODGreedy}
    \textbf{Input:} $N_\text{max} > 0$, finite set $\mathcal{P}_\text{train} \subset \mathcal{P},$ tolerance $\varepsilon$
	\begin{algorithmic}[1]
    \State Choose $\mu_1 \in \mathcal{P}_\text{train}$; set $\Phi_1 := \{\varphi_1^i\}_{i=1}^l = \text{POD} \left\{ U_h (\mu_1) \right\}$, $V_1:= \text{span}\langle \Phi_1 \rangle$; set $N:=1$
    \While{ $\max_{\mu \in \mathcal{P}_\text{train}} \Delta_e(V_N,\mu) > \varepsilon$ and $N < N_\text{max}$}
    \State Define $\mu_{N+1} := \text{argmax}_{\mu \in \mathcal{P}_\text{train}} \Delta_e(V_N,\mu)$ 
    \State Define $\{\varphi^i_{N+1}\}_{i=1}^l := \text{POD} \left\{  U_h (\mu_{N+1}) - \text{Proj}_{V_N} (U_h(\mu_{N+1} ))  \right\}$ 
    \State Define $\Phi_{N+1} := \text{orthonormalize} \left( \Phi_N \cup \{\varphi^i_{N+1}\}_{i=1}^l  \right)$
    \State Define $V_{N+1} := \text{span} \langle \Phi_{N+1} \rangle,$ Set $N = N+1$ 
    \EndWhile
    \State Define $V_r := V_N, N_r := \text{dim} (V_r)$ 
    \State \textbf{Return} $V_r, N_r;$
	\end{algorithmic} 
\end{algorithm}

\subsection{MOR for the wave model}
\label{sec:MORwaveModel}
In the following, we create a reduced order model (ROM) for the model of wave propagation in damaged material described in Section \ref{sec:waveModel} with fixed damage location, i.e.\ setting $x_\text{1D} = x_\text{2D} = 3$ to be constants and varying the damage properties $\mu = (c_\text{D}, d_\text{D})$.
We apply the POD-Greedy algorithm to construct a reduced basis which is orthonormal with respect to the weighted inner product
\begin{equation*}
    \langle v_1,v_2\rangle_M = v_1^T M v_2,
\end{equation*}
where $M$ is the mass matrix. 
We use the parameter set 
\begin{equation}
\label{eq:Ptrain}
    \mathcal{P}_\text{train} = \Big\{ (c_\text{D},d_\text{D}): c_\text{D} \in \{ 0.02,0.03,0.04,...,0.20 \}, d_\text{D} \in \{ 0.05,0.06,0.07,...,0.50 \} \Big\} 
\end{equation}
that contains $\# \mathcal{P}_\text{train} = 874$ parameters. The projection error at the end time $t=T,$ i.e.\ $k = m$ defined by
\begin{equation}\label{err:proj}
    e_\text{proj} := \sqrt{ (\textbf{u}^m_h(\mu ) - \textbf{u}^m_p (\mu ) )^T M (\textbf{u}^m_h (\mu ) - \textbf{u}^m_p (\mu ) ) }
\end{equation}
with $\textbf{u}_p^m(\mu)$ denoting the projection onto the current space $V_N$ serves as error indicator $\Delta_e(V_N,\mu)$ in the greedy step, line 3 of algorithm \ref{Alg:PODGreedy}. For a POD-Greedy approach with an error estimate that is independent of the full order dimension $N$ and based on an a-posteriori error estimation, we refer to \cite{Glas2020}. The error tolerance  is set to $\varepsilon=0.1$ and the maximum number of iterations is $N_\text{max}=1000$. In each iteration $l=1$ new basis vector is added. We perform weighted POD using the mass matrix as weighting matrix, compare \cite{Vol13}. \\
In practice, we first compute all high dimensional solutions $U_h(\mu)$ to the parameters $\mu \in \mathcal{P}_\text{train}$ and store them. This has the advantage that those high dimensional solutions do not need to be computed repeatedly in each iteration, this step took $1973  \text{s} \approx  33 \text{min}$. The POD-Greedy approach is completed after $1455 \text{s} \approx 24 \text{min}$  when the projection error on $\mathcal{P}_\text{train}$ reaches the error tolerance $\varepsilon = 0.1$. The computed basis of the RB space $V_r$ consists of 175 basis functions. We note that in order to decrease the computational effort in the offline phase  an error estimator according to \cite{Glas2020}, which is independent of the full order dimension, could be used.

\begin{figure}
\begin{minipage}{72mm}
\includegraphics[height=50mm]{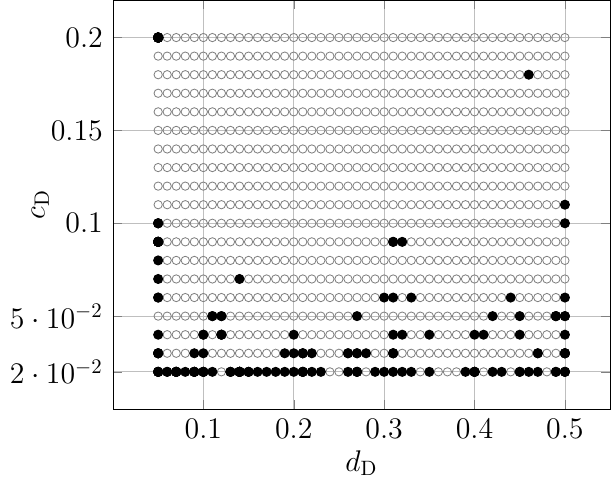}
\caption{Parameter set $\mathcal{P}_\text{train}$ used in experiment 1 represented by circles. The ticked circles represent parameters that were identified in the POD-Greedy algorithm to compute basis vectors, compare \cite[Abb.\ 3]{GH2024}. }
\label{fig:exp1:params}
\end{minipage}
\hfil
\begin{minipage}{65mm}
\includegraphics[height=50mm]{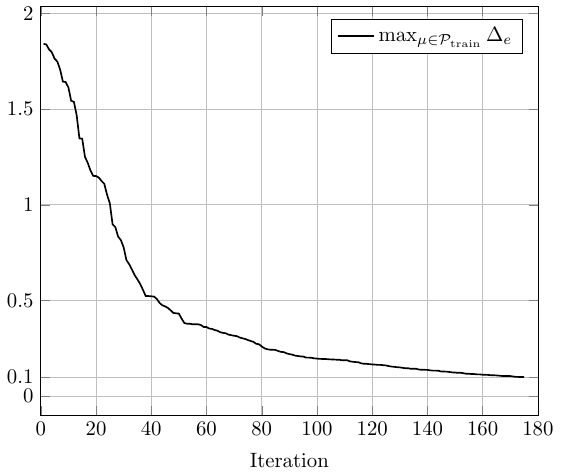}
\caption{Decay of the error indicator $\Delta_e(V_N,\mu) = e_\text{proj}$ evaluated over $\mathcal{P}_\text{train}$ within the POD-Greedy algorithm over the length of the reduced basis, which is equal to the iteration + 1, compare \cite[Abb.\ 3]{GH2024}.}
\label{fig:exp1:developmentProjectionError}
\end{minipage}
\end{figure}

\begin{figure}
\includegraphics[width=0.75\linewidth]{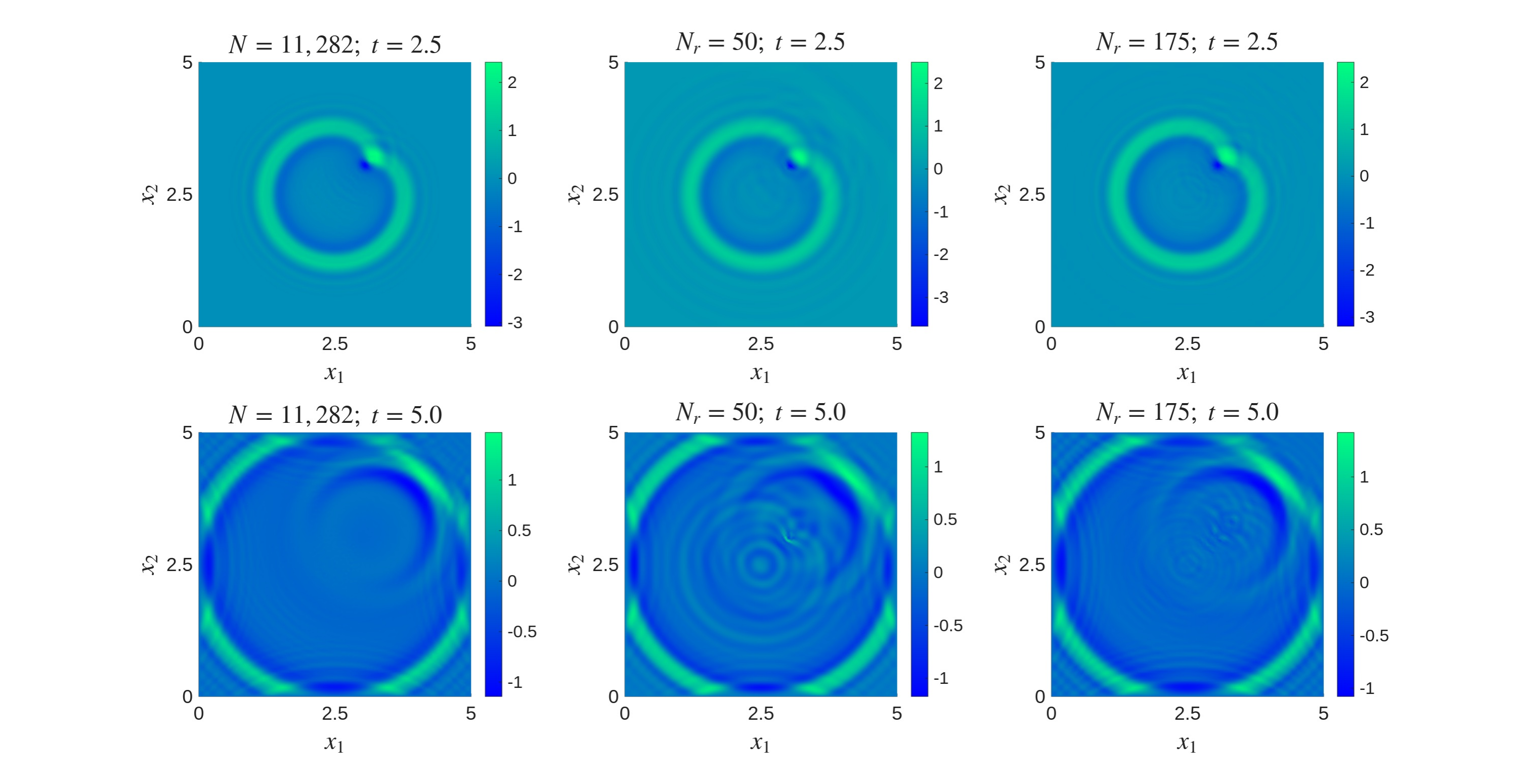}
\caption{Solution of the FOM for $\mu = (0.05,0.1)$ (left) and approximations computed with a ROM of dimension 50 (middle) and 175 (right) at times $t= 2.5$ (top) and $t = 5.0 $ (bottom).}
\label{fig:V1:differentROMs}
\end{figure}

In Fig.\ \ref{fig:exp1:params} the parameter values $\mu \in \mathcal{P}_\text{train}$ are represented as circles and ticked if the parameter was selected during the POD-Greedy algorithm, compare also \cite{GH2024}. Parameters with smaller values for $c_\text{D}$ were selected significantly more frequently than parameters with larger values.
This is due to the fact that if the propagation velocity in the damage is further reduced compared to the damage-free area, the damage has a stronger effect on the emitted wave resulting in dynamics that are more difficult to approximate. Out of the $\# \mathcal{P}_\text{train} = 874$ parameter values in the training set, 83 parameter values were selected and used to construct the reduced basis. \\
Figure~ \ref{fig:exp1:developmentProjectionError} shows the decay of the projection error in the course of the POD-Greedy algorithm, compare also \cite{GH2024}. The projection error falls rapidly in the first 40 iterations and decreases more slowly thereafter. An error bound $\varepsilon<0.1$ would presumably result in significantly more basis vectors and a considerably longer offline phase. \\
The solution to the full order model is illustrated in Figure~\ref{fig:V1:differentROMs} as well as its approximation computed by solving the reduced model for $N_r = 50$ and $N_r=175$ basis vectors and the parameter value $\mu = (0.05,0.1)$ at the time points $t = 2.5$ (top) and $t = 5$ (bottom). As expected, the approximation quality is better for $N_r = 175$. Nevertheless the essential structures are clearly visible even for $N_r=50$. \\
Table \ref{tab:Ch:V1:Zeitdauer} compares the high dimensional model and the reduced model in terms of the time required to compute the solutions for different parameter values.
The times needed to solve the high dimensional model and the reduced model are similar for all parameter values listed and takes 2.10s and 0.13s on average, respectively, making the reduced model 16 times faster than the high dimensional model.\\

In the next sections, parameter subsampling strategies are applied to the parameter set $\mathcal{P}_\text{train}$ \eqref{eq:Ptrain}. Based on the subsampled parameter set, different ROMs are set up and compared to the one from this section. In the following, we refer to the model reduction performed so far as \textit{experiment 1}.

\begin{table}
\caption{Computation times for the FOM and ROM for  different parameter values $\mu = (c_\text{D},d_\text{D})$  and the resulting acceleration factors.}
\label{tab:Ch:V1:Zeitdauer}
\begin{tabular}[h]{@{}lllll@{}}
\hline
$c_\text{D}$  & $d_\text{D}$  & Computation time &Computation time &Acceleration factor \\
 & & FOM in s & ROM in s & \\
\hline
0.02 & 0.05 & 2.07 &0.14& 14.8\\
0.02 & 0.25 & 2.05 & 0.13& 15.8\\
0.02 & 0.50 & 2.13 & 0.13 & 16.4\\
0.10 & 0.05 & 2.07 & 0.13 & 15.9 \\
0.10 & 0.25 & 2.12 & 0.14 & 15.1 \\
0.10 & 0.50 & 2.11 & 0.13 & 16.2 \\
0.20 & 0.05 & 2.09 & 0.13 & 16.1 \\
0.20 & 0.25 & 2.07 & 0.13 & 15.9 \\
0.20 & 0.50 & 2.18 & 0.12 & 18.2 \\
\hline
\end{tabular}
\end{table}

\section{Parameter Subsampling Strategies}

The duration of the offline phase and the quality of the ROM are highly influenced by the parameter set $\mathcal{P}_\text{train}.$ On the one hand, a finely sampled parameter set $\mathcal{P}_\text{train}$ ensures that parameters leading to a variety of dynamics are included and therefore, potentially improves the quality of the ROM. On the other hand, a large training set $\mathcal{P}_\text{train}$ prolongs the offline phase.   
\\
In Section \ref{sec:MORwaveModel}, 83 out of $\#\mathcal{P}_\text{train} = 874$ parameters were selected in line three of the POD-Greedy algorithm \ref{Alg:PODGreedy}. If we only use those 83 parameters to begin with, it takes 180s to compute all corresponding high dimensional solutions and 200s to complete the POD-Greedy algorithm. This leads to a significantly shorter offline phase of 380s $\approx$ 6.33min instead of 3425s $\approx$ 57.08min, while the quality of the ROM stays the same. This motivates the use of different subsampling strategies with the goal to identify a subset $\mathcal{P}_\text{k}$ of $\mathcal{P}_\text{train}$ that consists of parameters that are likely chosen in the POD-Greedy algorithm. In the following, different subsampling strategies are described and applied to create  $\mathcal{P}_k \subset \mathcal{P}_\text{train}.$  

\subsection{Subsampling based on DEIM and pivoted QR decomposition } 

The work \cite{Subsampling} introduces several subsampling strategies for the offline phase of an RB method based on Discrete Empirical Interpolation (DEIM), see algorithm \ref{alg:DEIM}, and the pivoted QR decomposition. The latter is a matrix decomposition of a matrix $\Tilde{\mathcal{U}}$ defined by
\begin{equation*}
    \Tilde{\mathcal{U}}^T \Pi = QR
\end{equation*}
consisting of an orthogonal matrix $Q$, an upper triangular matrix $R$ and a column permutation matrix $\Pi$. The pivots are given by the columns of the matrix $\Pi$. 
\\
Observations in the paper \cite{Subsampling} show similarities between the parameter values selected by the greedy algorithm and the pivots of the pivoted QR decomposition or the DEIM interpolation points. \\
\\
\begin{algorithm}
	\caption{DEIM, see \cite{Subsampling, CS2010}} 
    \label{alg:DEIM}
    \textbf{Input:} Snapshot matrix $\mathcal{F}$ , tolerance for the SVD $\varepsilon_\text{SVD}$ \\
    \textbf{Output:} $U,P,I$ 
	\begin{algorithmic}[1]
    \State Calculate the SVD of the matrix $\mathcal{F}$; yielding the first $l$ left singular vectors $\{ u_j\}_{j=1}^l$.
    \State Set $[\sim , p_1] = \max (|u_1 |).$
    \State Set $U=u_1, I = p_1$ and $P = [e_{p_1}].$
		\For {$i=2,\ldots, l$}
		      \State Solve $P^T U \alpha = P^T u_i.$ 
            \State Define the residual $r = u_i - U\alpha .$ 
            \State Set $[\sim , p_i] = \max (|r|)$ 
            \State Update $U:= [U u_i], I:= [I p_i]$ and $P:= [P e_{p_i}].$ 
		\EndFor
	\end{algorithmic} 
\end{algorithm}
The subsampling procedure consists of two phases. In the first phase, the POD-Greedy algorithm is started for a coarse error bound $\varepsilon_g$ on the fine training set $\mathcal{P}_\text{train}$. As soon as the coarse error bound $\varepsilon_g$ is reached, the second phase begins. The basis vectors computed so far are used to set up a first reduced model. Based on the reduced model, the snapshot matrix 
\begin{equation*}
    \Tilde{\mathcal{U}} = 
    \begin{pmatrix}
\textbf{u}^1_r (\mu_1)^T & ... &  \textbf{u}^m_r(\mu_1)^T\\
\vdots & \ddots & \vdots \\
\textbf{u}^1_r (\mu_p)^T & ... &  \textbf{u}^m_r(\mu_p)^T
\end{pmatrix}
\end{equation*}
is computed.
Then, 
\begin{enumerate}
    \item[a)] the pivoted QR decomposition is applied to the matrix $\Tilde{\mathcal{U}}^T$ and the first pivots are used to identify a subset $\mathcal{P}_k\subset \mathcal{P}_\text{train}$ or
    \item[b)] DEIM is applied to the matrix $\Tilde{\mathcal{U}}$ and the parameter values corresponding to the first interpolation points form $\mathcal{P}_k \subset \mathcal{P}_\text{train}.$
\end{enumerate}
The POD-Greedy algorithm is then continued on the parameter set $\mathcal{P}_k$ until the finer target tolerance $\varepsilon$ or the maximum number of iteration steps $N_\text{max}$ is reached.  The procedure is summarized in algorithm \ref{alg:Subsamplen}.
The final set of computed basis vectors is used to set up the ROM with the RB-Galerkin method described in Section \ref{subsec:GlaerkinRB}. In the following Section \ref{sec:ApplicationSubsampling} we refer to those strategies as \textit{QR, ROM} and \textit{DEIM, ROM}.

\begin{algorithm}
	\caption{Subsampling, compare algorithm 3.1 in \cite{Subsampling}} 
    \label{alg:Subsamplen}
    \textbf{Input:} Parameter set $\mathcal{P}_\text{train}$, tolerance $\varepsilon$ and coarse tolerance $\varepsilon_g$, maximum number of iterations $N_\text{max}$ \\
    \textbf{Output:} reduced basis $\Phi$
	\begin{algorithmic}[1]
    \State Choose $\mu_1 \in \mathcal{P}_\text{train}$ randomly and determine $\varphi_1 = \text{POD}\{ U_h (\mu_1) \}$; $\Phi_1=\{\varphi_1\},$ $V_1:= \text{span} \langle \Phi_1 \rangle$
    \Statex \textbf{Phase 1:} 
    \While {$\max_{\mu \in \mathcal{P}_\text{train}} \Delta_e(V_N,\mu) > \varepsilon_g$ and $l < N_\text{max}$}
    \State Define $\mu_{l+1} := \text{argmax}_{\mu \in \mathcal{P}_\text{train}} \Delta_e(V_N,\mu)$ 
    \State Define $\varphi_{l+1} := \text{POD} \left\{  U_h (\mu_{l+1}) - \text{Proj}_{V_l} (U_h(\mu_{l+1} ))  \right\}$ 
    \State Define $\Phi_{l+1} := \text{orthonormalize} \left( \Phi_l \cup \{ \varphi_{l+1} \} \right)$, $V_{l+1} := \text{span} \langle \Phi_{l+1} \rangle , l := l+1$
    \EndWhile  
    \Statex \textbf{Phase 2:}
    \State Calculate the pivoted QR decomposition of $\Tilde{\mathcal{U}}^T$ or apply DEIM to $\Tilde{\mathcal{U}}$. Identify the indices $I$ of the QR pivots or the DEIM interpolation points 
    \State Determine a new training set $\mathcal{P}_k$ using $I$ 
    \While {$\max_{\mu \in \mathcal{P}_k} \Delta_e(V_N,\mu)> \varepsilon$ and $l < N_\text{max}$}
    \State Define $\mu_{l+1} := \text{argmax}_{\mu \in \mathcal{P}_k} \Delta_e(V_N,\mu)$ 
    \State Define $\varphi_{l+1} := \text{POD} \left\{  U_h (\mu_{l+1}) - \text{Proj}_{V_l} (U_h(\mu_{l+1} ))  \right\}$ 
    \State Define $\Phi_{l+1} := \text{orthonormalize} \left( \Phi_l \cup \{ \varphi_{l+1} \} \right)$, $V_{l+1} := \text{span} \langle \Psi_{l+1}\rangle, l := l+1$ 
    \EndWhile 
    \State Define $\Phi:= \Phi_N$
    \State \textbf{Return} $\Phi_N$
	\end{algorithmic} 
\end{algorithm}

\noindent As a variant of the strategies described above, the FOM  is used to set up
\begin{equation*}
    \mathcal{U} = 
    \begin{pmatrix}
\textbf{u}^1_h (\mu_1)^T & ... &  \textbf{u}^m_h(\mu_1)^T\\
\vdots & \ddots & \vdots \\
\textbf{u}^1_h (\mu_p)^T & ... &  \textbf{u}^m_h(\mu_p)^T
\end{pmatrix} 
\end{equation*}
and the pivoted QR decomposition or the DEIM algorithm are applied to $\mathcal{U}$ instead of $\Tilde{\mathcal{U}}$ to subsample $\mathcal{P}_\text{train}$. 
We refer to those strategies as \textit{DEIM, FOM} and \textit{QR,FOM}.

\subsection{Polytope Division Method (PDM)}

\begin{figure}
\includegraphics[width=20mm,height=20mm]{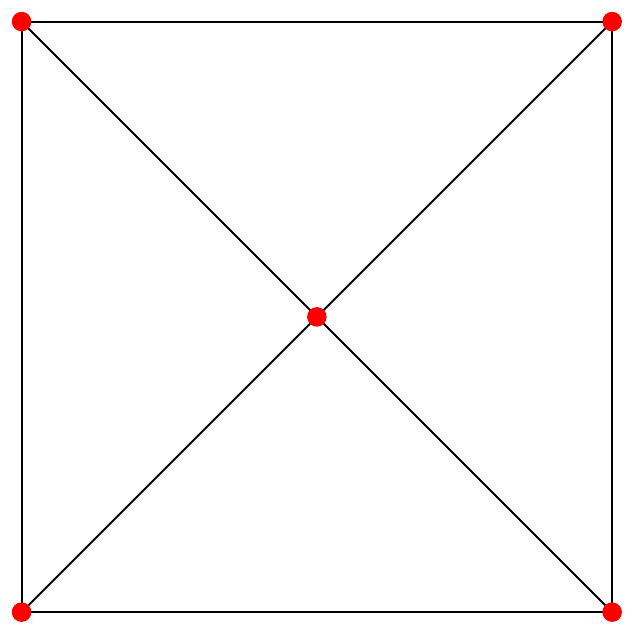}~1)
\hfil
\includegraphics[width=20mm,height=20mm]{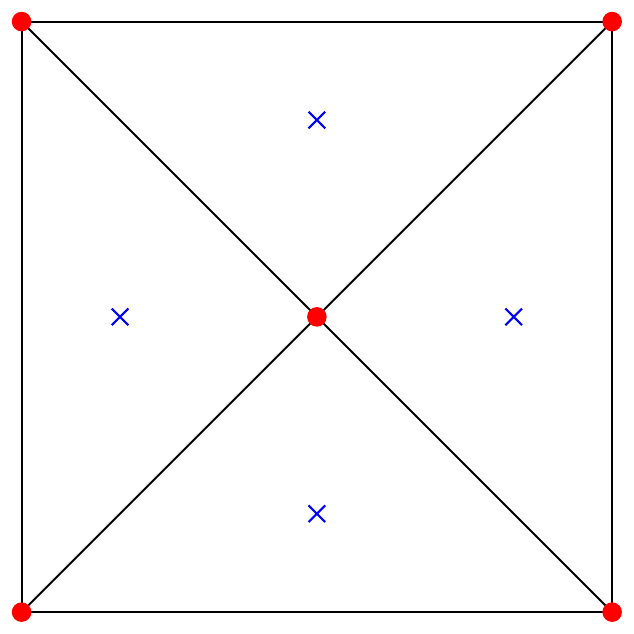}~2)
\hfil
\includegraphics[width=20mm,height=20mm]{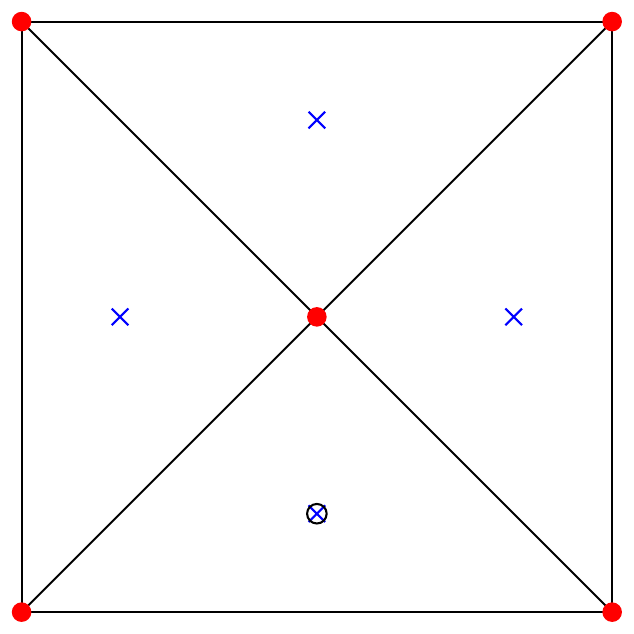}~3)
\hfil
\includegraphics[width=20mm,height=20mm]{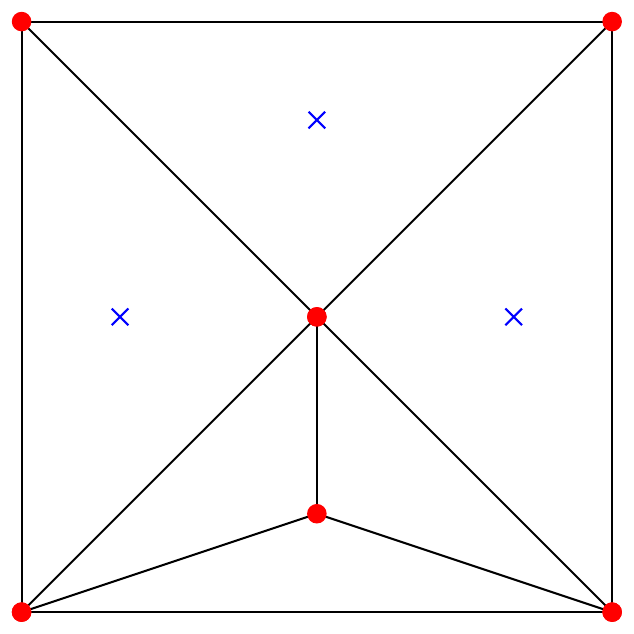}~4)
\hfil
\includegraphics[width=20mm,height=20mm]{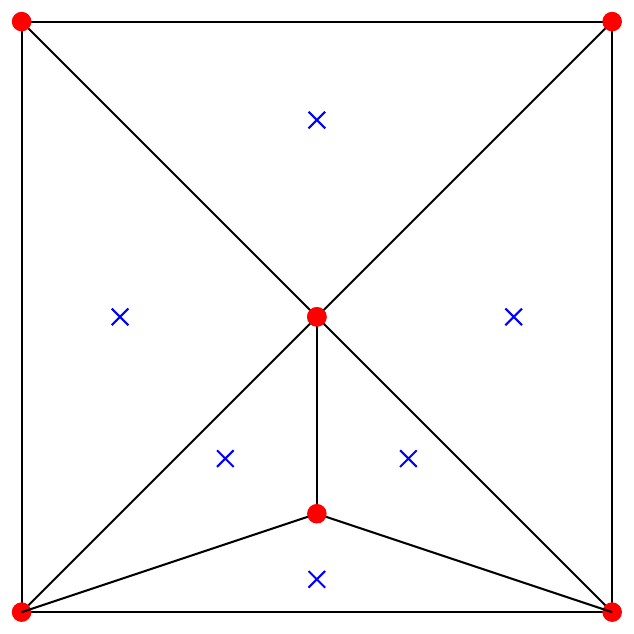}~5)
\hfil
\includegraphics[width=20mm,height=20mm]{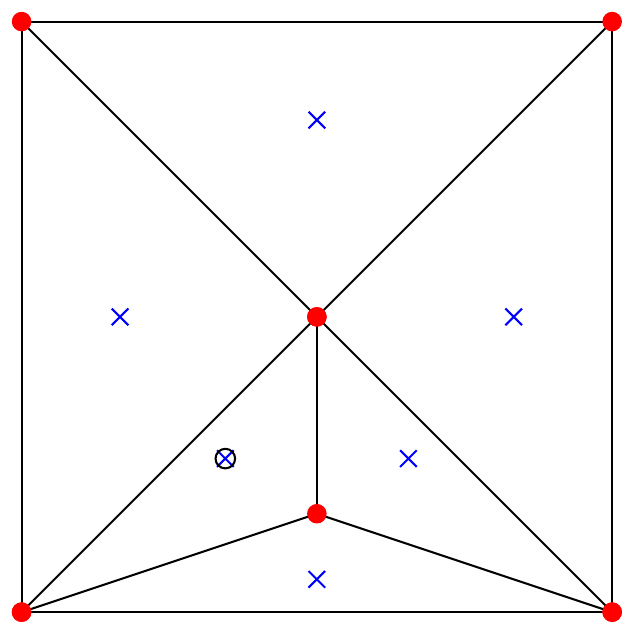}~6)
\caption{Illustration of steps 1-6 of the Polytope Division Method.}
\label{fig:PDM:idea}
\end{figure}

Another subsampling strategy is the Polytope Division Method (PDM) introduced in \cite{nielen2024polytope}. In strategy \textit{PDM} the parameter subset $\mathcal{P}_\text{k}$ is built iteratively as illustrated in Fig. \ref{fig:PDM:idea}. At first, the corners and the center of parameter set $\mathcal{P}$ form the so called set of nodes. Here, we follow the PDM approach using barycenters of the corresponding polytopes. A function to be minimized, in our case the projection error, is evaluated at the barycenters. The barycenter that maximizes the function is added to the sets of nodes. The barycenters of the resulting new polytopes are computed and the process is repeated  \cite{nielen2024polytope}.  \\

\section{Application and comparison of the subsampling strategies}
\label{sec:ApplicationSubsampling}

In \textit{experiment 1}, a total of 83 out of $\# \mathcal{P}_\text{train} = 874$ parameter values were used to set up a ROM of dimension 175, see Section \ref{sec:MORwaveModel}.
In the following, we apply the described subsampling strategies to the parameter set $\mathcal{P}_\text{train}$ and subsample sets $\mathcal{P}_k \subset \mathcal{P}_\text{train}$ that consist of 100 parameters and are used to set up a ROM of dimension 175 in each case. We evaluate the subsampling strategies based on three criteria: the duration of the offline phase, the subsampled set and the quality of the ROM as measured by the projection error. 

\subsection{Duration of the offline phase}
Of particular interest is whether the subsampling strategies shorten the offline phase, i.e.\ the time required to construct the ROM. Below, we detail the time needed for each step and the total offline phase for the various subsampling strategies.  \\ 
\\
We start the first phase of the strategies \textit{QR,ROM} and \textit{DEIM,ROM}, with computing and storing all high dimensional solutions to $\mathcal{P}_\text{train}$ once in advance which takes 1980s. Then, we start the POD-Greedy algorithm using the error indicator $\Delta_e(V_N,\mu) = e_\text{proj}$ and the coarse tolerance $\varepsilon_g = 1.$ The algorithm terminates after 25 iterations taking 160s.  \\
Phase 2 begins with setting up the matrices $\, \Tilde{\mathcal{U}}^T$ and $\, \Tilde{\mathcal{U}}$ in 80s. The calculation of the pivoted QR decomposition and the DEIM algorithm take 0.60s and 0.15s respectively. 
The parameter set $\mathcal{P}_\text{train}$ is subsampled based on the first 100 pivots in case of \textit{QR,ROM} or interpolation points in case of \textit{DEIM,ROM}.
We then continue the POD-Greedy algorithm on $\mathcal{P}_k$ until the desired basis length 175 is reached which takes 205s.  The total duration of strategies \textit{QR,ROM} and \textit{DEIM,ROM} is therefore
        \begin{align*}
             t_\text{QR,ROM} = t_\text{DEIM,ROM} = 1980s + 160s + 80s + 205s = 2425s,
        \end{align*}
neglecting the time required for the QR decomposition or for the DEIM algorithm. Strategy \textit{QR,FOM} involves setting up the matrix $\mathcal{U}^T$ and reducing it to the hundredth columns for capacity reasons, computing the QR decomposition of the matrix $\mathcal{U}^T_{100}$ and performing the POD-Greedy algorithm with the parameter values corresponding to the first 100 pivots to compute 175 basis functions. The time requirements are 2810s, 10s and 225s leading to a total of
        \begin{align*}
        t_\text{QR,FOM} = 2810s +10s + 225s = 3045 s.
        \end{align*}
Strategy \textit{DEIM,FOM} involves setting up the matrix $\mathcal{U}$, applying the DEIM algorithm to $\mathcal{U}_{100},$ and performing the POD-Greedy algorithm with the parameter values corresponding to the first 100 DEIM interpolation points. The individual steps take 2880s, 1s and 240s. The total duration, disregarding the duration of the DEIM algorithm, is therefore
    \begin{equation*}
        t_\text{DEIM,FOM} = 2880s+1s+240s = 3121s.
    \end{equation*}

For strategy \textit{PDM1}, we initialize the algorithm with the corners and the center of the parameter space $\mathcal{P}$ as nodes. Thus, the set $\mathcal{P}_k$ at first consists of the four corresponding barycenters. In each iteration of the POD-Greedy algorithm one new parameter is added to $\mathcal{P}_k$ following the PDM method until $\mathcal{P}_k$ contains 100 parameters. While doing so, we replace the computed barycenters by the closest value of $\mathcal{P}_\text{train}$ not already in $\mathcal{P}_k$, since we want to subsample $\mathcal{P}_\text{train}$ and not $\mathcal{P}$. The POD-Greedy algorithm stops as soon as 175 basis functions are computed, taking altogether
    \begin{equation*}
        t_\text{PDM1} = 553 s.
    \end{equation*}
For the \textit{PDM2} strategy , we use the final set $\mathcal{P}_k$ from strategy \textit{PDM1} right from the start of the POD-Greedy algorithm. It takes 384s to identify the 100 parameters following strategy \textit{PDM1} and an additional 220s to run the POD-Greedy algorithm on the subsampled parameter set $\mathcal{P}_k$, yielding a total of
    \begin{equation*}
        t_\text{PDM2} = 384s + 220s = 604s.
    \end{equation*}
Last, we consider a coarser uniform parameter set and a random parameter set of cardinality $\# \mathcal{P}_k = 100$ each, referred to as \textit{coarse} and \textit{random}. The offline phase took in both cases 
        \begin{equation*}
            t_\text{random} = t_\text{coarse} = 220s+ 215s = 415s.
        \end{equation*}
It takes 220s to compute the high dimensional solutions $U_h(\mu )$ for all 100 parameter values once in advance and an additional 215s to complete the subsequent POD-Greedy algorithm.\\ 
\\
All times are listed in Table \ref{tab:subsample:durationOnlinePhase}. We find that the strategies based on pivoted QR decomposition and DEIM shorten the  offline phase by 5-6min in case of \textit{QR, FOM} and \textit{DEIM, FOM} and by 20min in case of \textit{QR, ROM} and \textit{DEIM, ROM} compared to \textit{experiment 1}, where we used the fine parameter set $\mathcal{P}_\text{train}$ the whole time. Here, more than half of the total duration of the offline phase falls to computing the FOM solutions to $\mathcal{P}_\text{train}$. Since the storage capacities allow us to compute and store all high dimensional solutions, even for the finer set $\mathcal{P}_\text{train}$ once in advance, not that much of a speed up is generated by the use of a smaller set $\mathcal{P}_k$ later on.  \\
The strategies based on PDM, on the other hand, have significantly shortened the offline phase to around 10 min.

\subsection{Subsampled parameter sets}

\begin{figure}[hbt]
     \centering
     \begin{subfigure}[b]{0.32\textwidth}
         \centering         \includegraphics[width=0.85\textwidth]{Params_exp1.pdf}
         \caption{experiment 1}
     \end{subfigure}
    \hfill
     \begin{subfigure}[b]{0.32\textwidth}
         \centering        \includegraphics[width=0.85\textwidth]{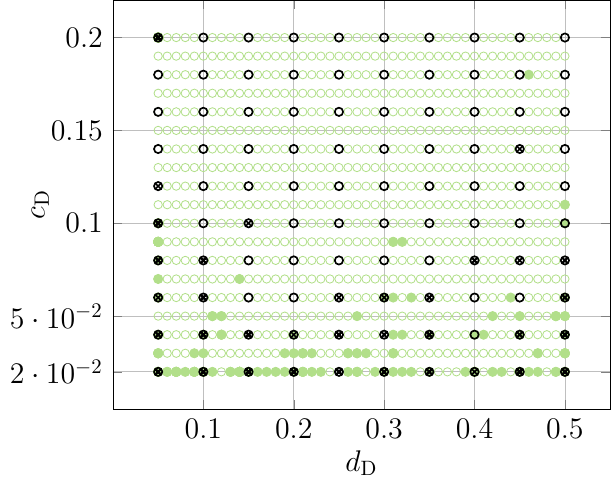}
         \caption{coarse}
     \end{subfigure}
     \begin{subfigure}[b]{0.32\textwidth}
         \centering        \includegraphics[width=0.85\textwidth]{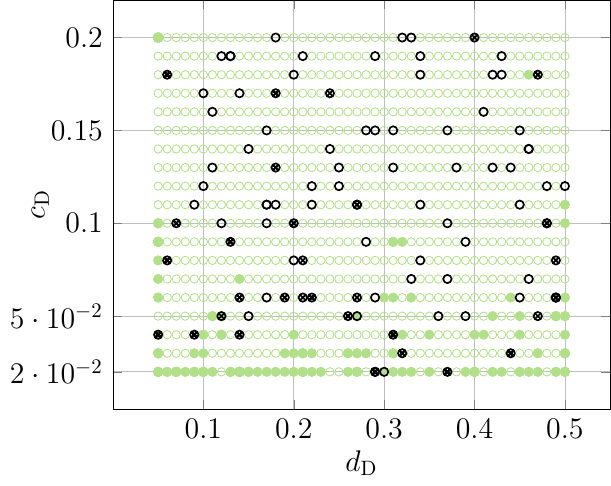}
         \caption{random}
     \end{subfigure}
     \hfill
     \begin{subfigure}[b]{0.32\textwidth}
         \centering  \includegraphics[width=0.85\textwidth]{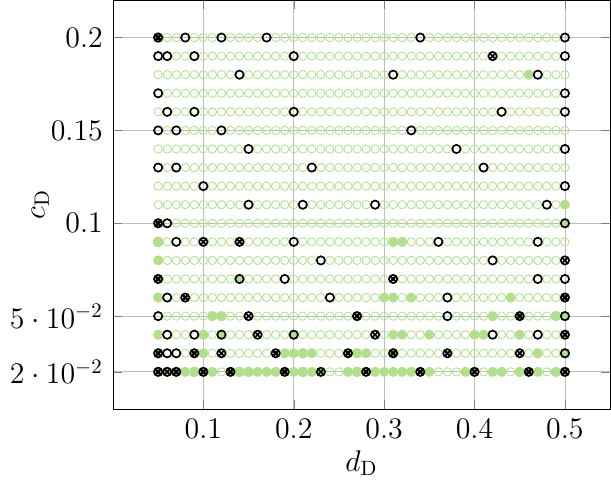}
         \caption{QR,ROM}
     \end{subfigure}
    \hfill
     \begin{subfigure}[b]{0.32\textwidth}
         \centering       \includegraphics[width=0.85\textwidth]{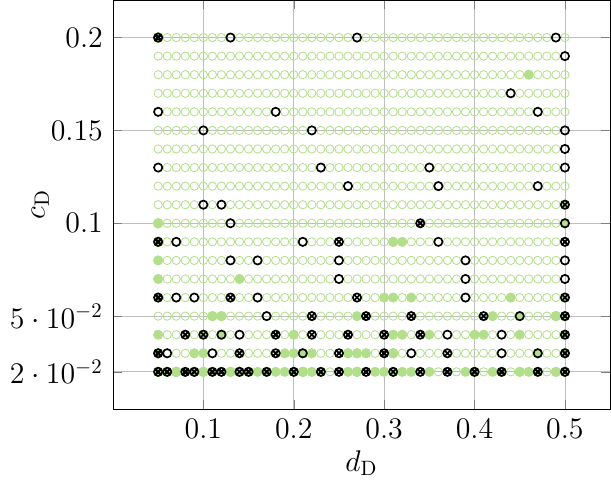}
         \caption{QR,FOM}
     \end{subfigure}
     \begin{subfigure}[b]{0.32\textwidth}
         \centering         \includegraphics[width=0.85\textwidth]{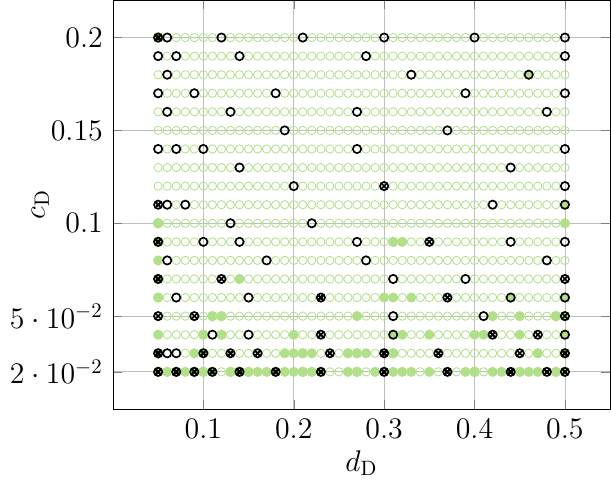}
         \caption{DEIM,ROM}
     \end{subfigure}
    \hfill
    \begin{subfigure}[b]{0.32\textwidth}
         \centering         \includegraphics[width=0.85\textwidth]{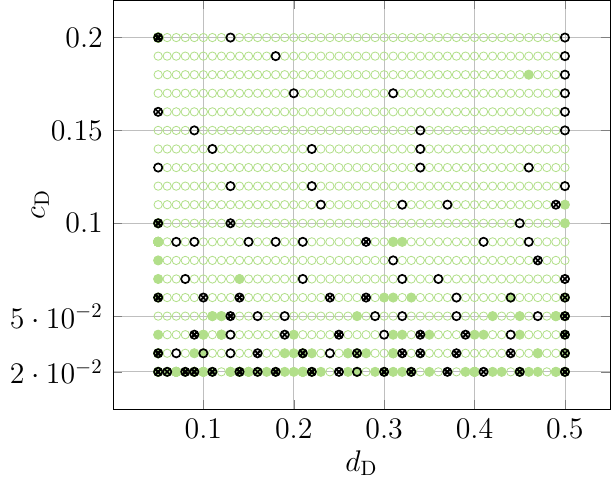}
         \caption{DEIM,FOM}
     \end{subfigure}
     \hfill
     \begin{subfigure}[b]{0.31\textwidth} 
        \hspace*{0.65cm}\includegraphics[width=0.85\textwidth]{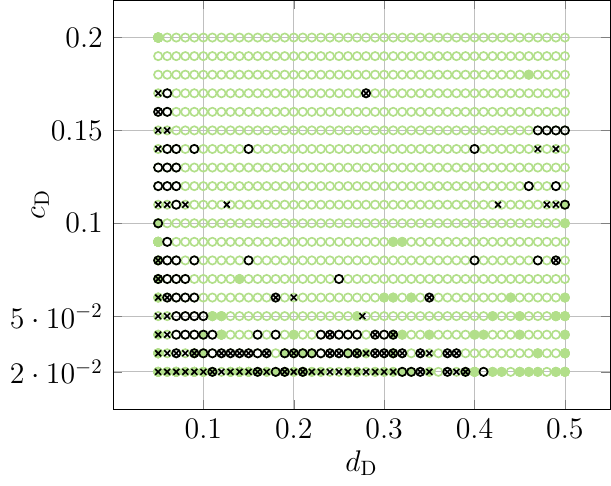}
         \caption{PDM 1}
     \end{subfigure}
    \hfill
     \begin{subfigure}[b]{0.305\textwidth}
         \hspace*{0.3cm}        \includegraphics[width=0.85\textwidth]{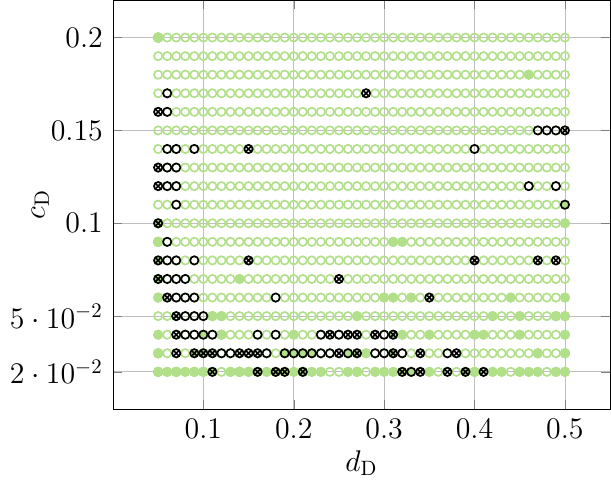}
         \caption{PDM 2}
     \end{subfigure}
        \caption{Subsampled parameter sets $\mathcal{P}_k \subset \mathcal{P}_{\text{train}}$ for the different strategies.}\label{fig:subsample:params}
\end{figure}

Figure~ \ref{fig:subsample:params} illustrates the different parameter sets used in the POD-Greedy algorithm. The subsampled parameter sets $\mathcal{P}_\text{k}$ are marked with black circles that are ticked, if the corresponding parameter has been chosen at least once in the POD-Greedy algorithm.
The green circles represent the full parameter set $\mathcal{P}_\text{train},$ the circles of the parameters chosen in \textit{experiment 1} are filled. \\
The results show that the greedy algorithm predominantly selected parameter values with small wave propagation velocity $c_\text{D}$. This is consistent with the results from Section \ref{sec:MORwaveModel} for \textit{experiment 1}. As mentioned before, a further reduction of the wave propagation velocity at the center of the damage leads to more evolved dynamics that are harder to approximate. Subsequently, the corresponding parameter values are more likely picked in the POD-Greedy algorithm. The PDM strategies were best able to identify the trend that parameter values $\mu \in \mathcal{P}_\text{train}$ with smaller $c_\text{D}$ are more important regarding the POD-Greedy algorithm. Also the strategies \textit{DEIM,FOM} and \textit{QR,FOM} predominantly subsampled parameters corresponding to smaller $c_\text{D}$ values. In the case of the strategies \textit{DEIM,ROM} and \textit{QR,ROM} this tendency is hardly recognizable and the selected parameter values are more equally distributed in the parameter domain. The first basis vectors primarily describe the emitted wave. The solutions of the ROM built on those basis vectors apparently do not show enough parameter specific differences to subsample the parameter sets.

\subsection{Error evaluation}
\begin{figure}
\begin{minipage}{95mm}
    \includegraphics[scale=0.6]{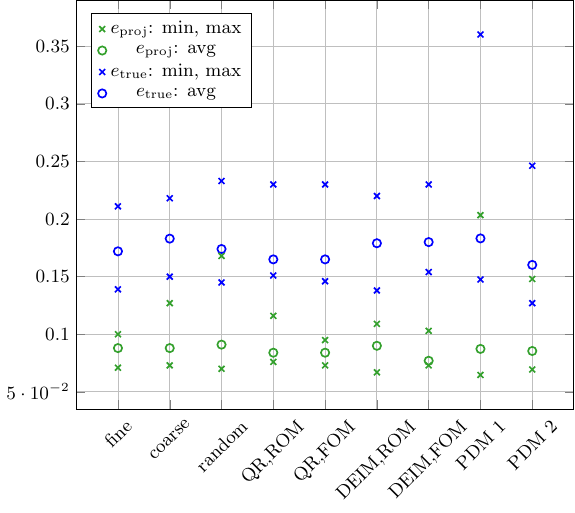}
    \caption{Comparison of the errors $e_\text{proj}$ and $e_\text{true}$ evaluated on the parameter set $\mathcal{P}_\text{train}$ over the subsampling strategies used to create the corresponding ROM. }
    \label{fig:errorComparison}
\end{minipage}\hspace{0.5cm}
%\hfill
\begin{minipage}{60mm}
%\setfloattype{table}
\caption{Duration of the offline phases for the different subsampling strategies.}
\label{tab:subsample:durationOnlinePhase}
\begin{tabular}{@{}ll@{}}
\hline
Strategy & Duration offline \\
& phase in min \\
\hline
            experiment 1 & 57 \\
            QR, ROM & 40.41 \\
            DEIM, ROM & 40.41 \\
            QR, FOM & 50.75 \\
            DEIM, FOM & 52.00 \\
            PDM 1 & 9.22 \\
            PDM 2 & 10.07 \\
            coarse & 6.92 \\
            random & 6.92 \\
\hline
\end{tabular}
\end{minipage}
\end{figure}

We evaluate the projection error $e_\text{proj}$ \eqref{err:proj} and the true error
\begin{equation}\label{err:true}
    e_\text{true} := \sqrt{ (\textbf{u}^m_h(\mu ) - \textbf{u}^m_a (\mu ) )^T M (\textbf{u}^m_h (\mu ) - \textbf{u}^m_a (\mu ) ) }
\end{equation}
at the end time $t=T,$ i.e.\ $k = m$ on the parameter set $\mathcal{P}_\text{train}.$ The resulting absolute minimum, 
$
    \min_{\mu \in \mathcal{P}_\text{train}} e_\text{proj/true},
$
and maximum, 
$
    \max_{\mu \in \mathcal{P}_\text{train}} e_\text{proj/true},
$
as well as the arithmetic mean
        \begin{equation*}
            \frac{1}{\# \mathcal{P}_\text{train}}\sum_{\mu \in \mathcal{P}_\text{train}} e_\text{proj/true}
        \end{equation*}
are visualized in Fig. \ref{fig:errorComparison}. Out of the tested subsampling strategies, the strategy \textit{PDM 1 } performed worst. A reason for that might be that the training set was computed iteratively beginning with only four parameters. Thus, the set was significantly smaller than the other sets at the start giving the POD-Greedy algorithm fewer parameters to choose from. Overall the strategies show a similar performance regarding the projection error and the true error. 
In case of strategy $\textit{QR,FOM}$ the tolerance $\varepsilon = 0.1$ is reached on all $\mathcal{P}_\text{train}$ even though the whole POD-Greedy algorithm is based on a subset $\mathcal{P}_k\subset \mathcal{P}_\text{train}.$ It is also the only strategy where the average true error falls below the value for \textit{experiment 1}.

\section{Conclusions and related work }

We considered a parameterized wave equation where a damage was modeled as a local reduction of the wave propagation velocity. The damage is characterized by two parameters, $c_\text{D}$ and $d_\text{D}$. We performed model order reduction for the wave model, considering a fixed damage position $(x_\text{1D},x_\text{2D}) = (3,3)$ and variable damage wave propagation $c_\text{D}$ and size $d_\text{D}$. We applied a reduced basis method with POD-Greedy algorithm and reduced the dimension from 11,282 to 175, leading to an acceleration factor of around 16. \\
Then, we applied various subsample strategies that, despite major differences in terms of duration, delivered similar results regarding the projection error and the true error on $\mathcal{P}_\text{train}$.  The parameter sets identified with strategies \textit{QR,FOM} and \textit{DEIM, FOM} contain many parameters with smaller values for $c_D.$ The latter lead to less well approximable high dimensional solutions, making strategies \textit{QR,FOM} and \textit{DEIM, FOM} suitable for identifying parameter values that lead to dynamics that are difficult to approximate. Such a tendency was not clearly observed for strategies \textit{QR,ROM} and \textit{DEIM, ROM} based on solutions of a reduced model. The strategies \textit{random} and \textit{coarse} led to similarly good results as the more evolved strategies with a significantly shorter offline phase. \\
\\
Finally, we want to refer the interested reader to related work done by the authors on the topic of MOR in the context of damage identification. 
In \cite{MasterthesisHannah}, additionally to the results presented here, MOR for the case of variable damage position $(x_\text{1D}, x_\text{2D})$ and constant $(c_\text{D}, \, d_\text{D})$ is examined and different strategies to perform the POD-Greedy algorithm stepwise are explored.  The work \cite{GH2024} gives an overview on the application of data assimilation methods for parameter estimation in the context of damage identification. It also presents the application of data-driven reduction methods to the numerical simulation of guided ultrasonic wave propagation in fiber metal laminates with damage. In addition, the work \cite{MGRMLL23} studies Bayesian analysis to infer damage parameters. In \cite{HFBRG25} apart from the wave model presented here, two additional models are considered which are a damaged fiber metal laminate model where the equation of linear elasticity is solved numerically and a second model which is a nonlinear elastic model of Aluminum plate with a damage. Both intrusive and non-intrusive model reduction approaches are applied to a fixed damage configuration.

\begin{acknowledgement}
 The authors acknowledge funding by the German Research Foundation (Deutsche Forschungsgesellschaft, DFG), project number 418311604 (Research Unit FOR3022: Ultrasonic Monitoring of Fibre Metal Laminates Using Integrated Sensors).
\end{acknowledgement}

\end{document}